\documentclass[11pt]{article}

\usepackage[margin=1in]{geometry}
\usepackage{amsmath,amssymb,amsthm,mathtools}
\usepackage{xcolor}
\usepackage[colorlinks=true,linkcolor=blue,citecolor=blue,urlcolor=blue,
  pdftitle={A Proof of Stein's Equidissection Conjecture},
  pdfauthor={Dmitry Krachun}]{hyperref}
\usepackage{microtype}

\usepackage{subcaption}
\usepackage{tikz}
\usetikzlibrary{arrows.meta}

\definecolor{markcol}{HTML}{087E8B}
\definecolor{edgecol}{HTML}{4B5563}
\definecolor{branchfill}{HTML}{CDE9EA}

\numberwithin{equation}{section}

\theoremstyle{plain}
\newtheorem{theorem}{Theorem}

\newtheorem{proposition}[theorem]{Proposition}

\newtheorem{conjecture}[theorem]{Conjecture}

\theoremstyle{definition}

\newtheorem{remark}[theorem]{Remark}
\newtheorem{example}[theorem]{Example}

\title{A Proof of Stein's Equidissection Conjecture}
\author{Dmitry Krachun%
  \thanks{Princeton University. Email:
  \href{mailto:dk9781@princeton.edu}{\texttt{dk9781@princeton.edu}}}}
\date{}

\begin{document}

\maketitle
\vspace{-1.5em}

\begin{abstract}
In 1970 Monsky showed that a square cannot be dissected into an odd number of equal-area triangles. His proof used an ingenious colouring of the plane based on a 2-adic valuation on $\mathbb{R}$ together with Sperner's lemma. Subsequent work established similar results for regular polygons, centrally symmetric polygons, and polyominoes. In 2000, Stein proposed a common generalisation of the centrally symmetric and polyomino cases to balanced polygons. In 2014, Rudenko introduced a new perspective on the problem and proved Stein's conjecture for dissections with rational vertex coordinates. Building on and streamlining his approach, we resolve the conjecture in full.

\end{abstract}

\section{Introduction}

An \emph{equidissection} of a polygon is a dissection into triangles of equal area. An \emph{odd equidissection} is an equidissection into an odd number of triangles. The basic question is then whether a given polygon admits such a dissection. 

Despite its simple formulation, the question of whether a square admits an odd equidissection was apparently not discussed until 1965, and only appeared in print in 1967 when Richman and Thomas asked it in the \emph{American Mathematical Monthly} \cite{RichmanThomas1967}. A year later, Thomas \cite{Thomas1968} proved that no such equidissection of the unit square exists if all vertices have rational coordinates with odd denominators. In 1970 Monsky famously solved the full problem by cleverly combining a colouring of the plane in three colours based on a $2$-adic valuation on $\mathbb{R}$ with Sperner's lemma \cite{Monsky1970}.

The later history broadened both the class of polygons and the available methods. In 1989 Kasimatis showed that any equidissection of a regular $n$-gon with $n\geq 5$ must consist of a number of triangles divisible by $n$ \cite{Kasimatis1989}. In the same year, Stein conjectured that no centrally symmetric $n$-gon can admit an odd equidissection and proved this assertion for $n=6$ and $n=8$ \cite{Stein1989Octagons}. A year later, Monsky proved Stein's conjecture for centrally symmetric polygons \cite{Monsky1990}; that proof is considerably more involved than the square proof, combining the three-colouring with a nontrivial homological argument.

In 1999 Stein proved that no polyomino of odd area admits an odd equidissection \cite{Stein1999Polyomino}, and in 2002 Praton extended the result to all polyominoes \cite{Praton2002}. For historical background and a survey of earlier results, see \cite{Stein2004Survey}.

A polygon is \emph{balanced}\footnote{Interestingly, balanced polygons also appear in the Hadwiger--Glur theorem \cite{HadwigerGlur1951}, which shows, as a special case, that balanced polygons are precisely those which are equidecomposable by translations with a square of the same area.} if, for every direction, the sum of the oriented boundary vectors parallel to that direction is zero. Equivalently, after subdividing boundary edges if necessary, its boundary segments can be paired into opposite vectors. This is Stein's notion of a special polygon \cite{Stein2000GeneralizedConjecture}. In 2000, Stein formulated the following general conjecture \cite{Stein2000GeneralizedConjecture}:

\begin{conjecture}[Stein]\label{conj:stein}
A balanced polygon admits no equidissection into an odd number of triangles.
\end{conjecture}

In particular, this conjecture includes all centrally symmetric polygons and all polyominoes as special cases. In 2012, Rudenko proved the special case of this conjecture for balanced lattice polygons of odd area \cite{RudenkoBalancedPolygons} and two years later for any balanced polygon for which the vertices of the equidissection (including the vertices of the polygon) are rational \cite{RudenkoTrivalentGraphs} (we cite the 2025 revision, which substantially streamlines the proof and exposition). The latter work introduced a simple reduction to a question about $\mathbb{Q}^2$-valued flows on trivalent graphs.

The main result of this paper is the full resolution of Stein's conjecture. As a consequence, we obtain a short unified proof for centrally symmetric polygons and all polyominoes.  

\begin{theorem}\label{thm:main}
    Stein's conjecture holds: no balanced polygon in the plane admits an odd equidissection.
\end{theorem}

The proof follows Rudenko's reduction from dissections into triangles to flows on trivalent graphs \cite{RudenkoTrivalentGraphs}, but the rest of the argument is substantially different. When all coordinates of the equidissection are rational, Rudenko works over the localisation $\mathbb Z_{(2)}$, whose residue field is $\mathbb{F}_2$; for a general number field $K$ generated by the coordinates of the equidissection, we work over the localisation $\mathcal{O}_{K, \mathfrak{p}}$ of $\mathcal{O}_K$ at a prime ideal $\mathfrak{p}$ above $2$, whose residue field may be larger. The crucial property of the residue field underlying the cycle decomposition of \cite{RudenkoTrivalentGraphs} is that three nonzero elements of $\mathbb{F}_2$ cannot sum to zero. Failure of this property for larger residue fields then forces us to develop a more general cancellation theorem whose proof uses the geometry of rank-two lattices over the localisation $\mathcal{O}_{K, \mathfrak{p}}$ together with a homological lifting argument.

The rest of the article is organised as follows. In Section~\ref{sec:flow-theorem} we state and prove a theorem about $K^2$-valued flows on a trivalent graph, where $K$ is a number field. We then prove our main theorem in Section~\ref{sec:proof-of-the-main-theorem}.

\medskip
\noindent \textbf{Acknowledgements:} I would like to thank Daniil Rudenko for introducing me to the problem of equidissection and for many enlightening discussions. I would also like to thank Alexey Balitskiy, Stepan Nesterov, Fedor Petrov, and Mura Yakerson for various comments on earlier versions of the paper. Additionally, I thank Micha\l{} Szachniewicz for his useful comments on the Tarski--Seidenberg theorem. 

\medskip
\noindent\textbf{Use of AI:}
The proof of the main theorem was obtained by the author using human intelligence. AI tools were used to assist with figure preparation and revising the exposition.

\section{Flows on trivalent graphs}\label{sec:flow-theorem}

Throughout this section, a graph means a finite loopless multigraph. The main result of this section is a generalisation of \cite[Theorem 1]{RudenkoTrivalentGraphs} to the case of arbitrary number fields. Recall the notion of a flow. For a trivalent graph $\Gamma=(V, E)$, let $\overrightarrow{E}$ be the set of its oriented edges, that is, pairs consisting of an edge and its orientation. A flow on $\Gamma$ taking values in an abelian group $A$ is defined to be a function $f:\overrightarrow{E}\rightarrow A$ satisfying 
\begin{enumerate}
    \item[(i)] $f(e)+f(e')=0$ for every pair of opposite oriented edges
    $e,e'\in \overrightarrow{E}$;

    \item[(ii)] $f(e_1)+f(e_2)+f(e_3)=0$ for every triple of distinct oriented edges
    $e_1,e_2,e_3\in \overrightarrow{E}$ with common source vertex.
\end{enumerate}

Let $K$ be a number field, let $\mathfrak{p}$ be a prime ideal in $\mathcal{O}_K$ dividing 2 and let $\nu_{\mathfrak{p}}:K\rightarrow\mathbb{Z}\cup\{\infty\}$ be the corresponding discrete valuation extending, up to scaling, the usual 2-adic valuation on $\mathbb{Q}$. Consider a trivalent graph $\Gamma$ with a flow $f$ valued in $A= K^2$. Let $\omega\in(\Lambda^2 K^2)^\vee$ be the standard area form given by
\[
\omega((x_1,y_1),(x_2,y_2))=x_1y_2-x_2y_1.
\]
For a vertex $v\in V$ define the number $A_v:=\omega(f(e_1),f(e_2))$, where $e_1,e_2$ are any two of the three oriented edges with source $v$. This quantity is well defined up to sign by the antisymmetry of $\omega$ and the flow condition. We make an arbitrary choice of the sign as it does not impact any of the claims. Indeed, since $\nu_{\mathfrak p}(2)\ge 1$, changing the sign of $A_v$ changes the sum in Theorem~\ref{thm:flow-main} by an element of strictly larger valuation.

The weight of a vertex $v\in V$ is then defined as $m_v:=\nu_{\mathfrak{p}}(A_v)$. The weight of the graph $\Gamma$ with the flow $f$ is the minimum of the vertex weights:
\[
m(\Gamma,f)=\min_{v\in V} m_v.
\]

The main result of this section is the following.

\begin{theorem}\label{thm:flow-main}
Let $\Gamma$ be a trivalent graph, let $K$ be a number field, and let $\mathfrak{p}$ be a prime ideal in $\mathcal{O}_K$ dividing 2. Given a flow
$f:\overrightarrow{E}\rightarrow K^2$ such that $m(\Gamma,f)<\infty$, one has 
\[
\nu_{\mathfrak{p}} \left(\sum_{v\in V\,:\, m_v=m(\Gamma, f)} A_v\right)>m(\Gamma, f).
\]
\end{theorem}
\begin{remark}
    Note that in the case $K=\mathbb{Q}$ the statement is equivalent to the assertion that the number of vertices $v\in V$ such that $m_v=m(\Gamma,f)$ is even, recovering \cite[Theorem 1]{RudenkoTrivalentGraphs}. In the case of a general number field, however, parity of the number of vertices $v$ satisfying $m_v=m(\Gamma, f)$ is no longer the relevant quantity.
\end{remark}

Before proving Theorem~\ref{thm:flow-main}, we need some preliminaries. Let $\mathcal{O}=\mathcal O_{K,\mathfrak p}$ be the localisation of the ring of integers of $K$ at a prime $\mathfrak p\mid 2$, let $\mathfrak m$ be the maximal ideal of $\mathcal{O}$, choose a uniformiser $\pi$ (i.e. an element with $\nu_{\mathfrak{p}}(\pi)=1$), and let $k=\mathcal{O}/\mathfrak m$ be the residue field. We work with rank-two $\mathcal{O}$-lattices in $K^2$. If $L$ is such a lattice, choose any $\mathcal{O}$-basis $a,b$ of $L$. The principal fractional ideal generated by $\omega(a,b)$ is independent of the basis; denote it by $\det L$. We write $\nu_{\mathfrak{p}}(L)$ for its valuation, equivalently $\nu_{\mathfrak{p}}(L)=\nu_{\mathfrak{p}}(\omega(a,b))$.

In the following, we fix the graph $\Gamma$, the number field $K$ and a prime $\mathfrak p$ lying over 2. We also fix an ambient lattice $P$ containing all vectors $f(e)$ with $e\in \overrightarrow{E}$. For each $v\in V$ let $L_v$ be the $\mathcal{O}$-submodule spanned by $f(e_1), f(e_2), f(e_3)$ where $e_1, e_2, e_3$ are the three oriented edges with source $v$; by the flow condition this $\mathcal{O}$-submodule is also spanned by any two of these vectors. If $A_v\neq 0$, then $L_v$ is a rank-two lattice, $\det L_v=(A_v)$, and consequently $\nu_{\mathfrak p}(L_v)=m_v$; if $A_v=0$, then $L_v$ has rank at most one.

\begin{remark}
Though we never explicitly use it in the proof, it is useful to keep in mind the structure of rank-two $\mathcal{O}$-lattices in $K^2$ contained in a given lattice $P$. Specifically, such lattices may be visualised using the Bruhat--Tits tree\footnote{Strictly speaking, the Bruhat--Tits tree is typically defined after completing $K$ at $\mathfrak p$. We work instead over $\mathcal O=\mathcal O_{K,\mathfrak p}$, but obtain the same tree, since its finite balls are determined by the quotients $\mathcal{O}/\pi^n\mathcal{O}$, which are unchanged upon completion.} of $\mathrm{PGL}_2(K_{\mathfrak p})$. Its vertices are homothety classes of rank-two $\mathcal O$-lattices, and two vertices are adjacent if they have representatives $L,L'$ satisfying $\pi L\subsetneq L'\subsetneq L$. After fixing a lattice $P$, every homothety class has a unique representative $L\subseteq P$ with $L\nsubseteq\pi P$, which we call its normalised representative; a lattice of this form is called normalised (with respect to $P$). The tree is $(|k|+1)$-regular. Given a vertex of the tree at distance $\ell$ from the vertex corresponding to the lattice $P$, its normalised representative $L$ has determinant $\pi^\ell\cdot \det P$. Moreover, for every vector $a\in P\setminus \pi P$, the set of all lattices $L$, again with $L\subseteq P$ and $L\nsubseteq\pi P$, containing $a$ is given by an infinite geodesic ray on the tree starting at $[P]$. We refer to \cite[Chapter~II, Section~1]{SerreTrees} for the proofs and more details. 
\end{remark}

\begin{example}
  Take $K=\mathbb{Q}(\sqrt5)$. The prime $2$ is inert in $K$, so we may take $\mathfrak p=(2)$ and $\pi=2$. With $\varphi=(1+\sqrt5)/2$ and $\alpha$ the image of $\varphi$ in the residue field $k$, one has
\[
k=\mathbb F_2(\alpha)\cong\mathbb F_4,
\qquad \alpha^2+\alpha+1=0,
\]
so the tree is $5$-regular. For $P:=\mathcal{O} a\oplus \mathcal{O} b$, there are exactly five normalised lattices corresponding to the five neighbouring vertices of $[P]$. Specifically, these are $\langle 2a, b\rangle, \langle a, 2b\rangle, \langle a+b, 2b\rangle, \langle a+\varphi b, 2b\rangle$, and $\langle a+(1+\varphi)b, 2b\rangle$.  
\end{example}

The proof of Theorem~\ref{thm:flow-main} is by induction on $m(\Gamma, f)$ and relies on the following two propositions: the first handles the base case and the main part of the induction step, while the second supplies the homological lifting that completes the step.

\begin{proposition}\label{prop:base-case}
Let $L \subseteq P$ be a lattice with $L \nsubseteq \pi P$ and assume that no vertex $v\in V$ has $L_v\supsetneq L$. Let $S_L=\{v \in V : L_v=L\}$ and define $\Gamma_L$ to be the subgraph of $\Gamma$ spanned by the edges $e\in E$ satisfying\footnote{Since $L$ and $\pi P$ are stable under negation, the condition $f(e)\in L\setminus \pi P$ does not depend on the orientation of $e$; the same applies to the similar conditions below.} $f(e)\in L\setminus \pi P$. Then the following statements hold:
\begin{enumerate}
    \item[(i)] The graph $\Gamma_L$ has vertices of degree 2 and 3 only;
    \item[(ii)] For every edge $e\in E(\Gamma)\setminus E(\Gamma_L)$ adjacent to a vertex in $V(\Gamma_L)$ one has $f(e)\in L\cap \pi P$;
    \item[(iii)] One has $S_L\subseteq V(\Gamma_L)$ and $\sum_{v\in S_L} A_v \in \pi \det L$.
\end{enumerate}

\end{proposition}

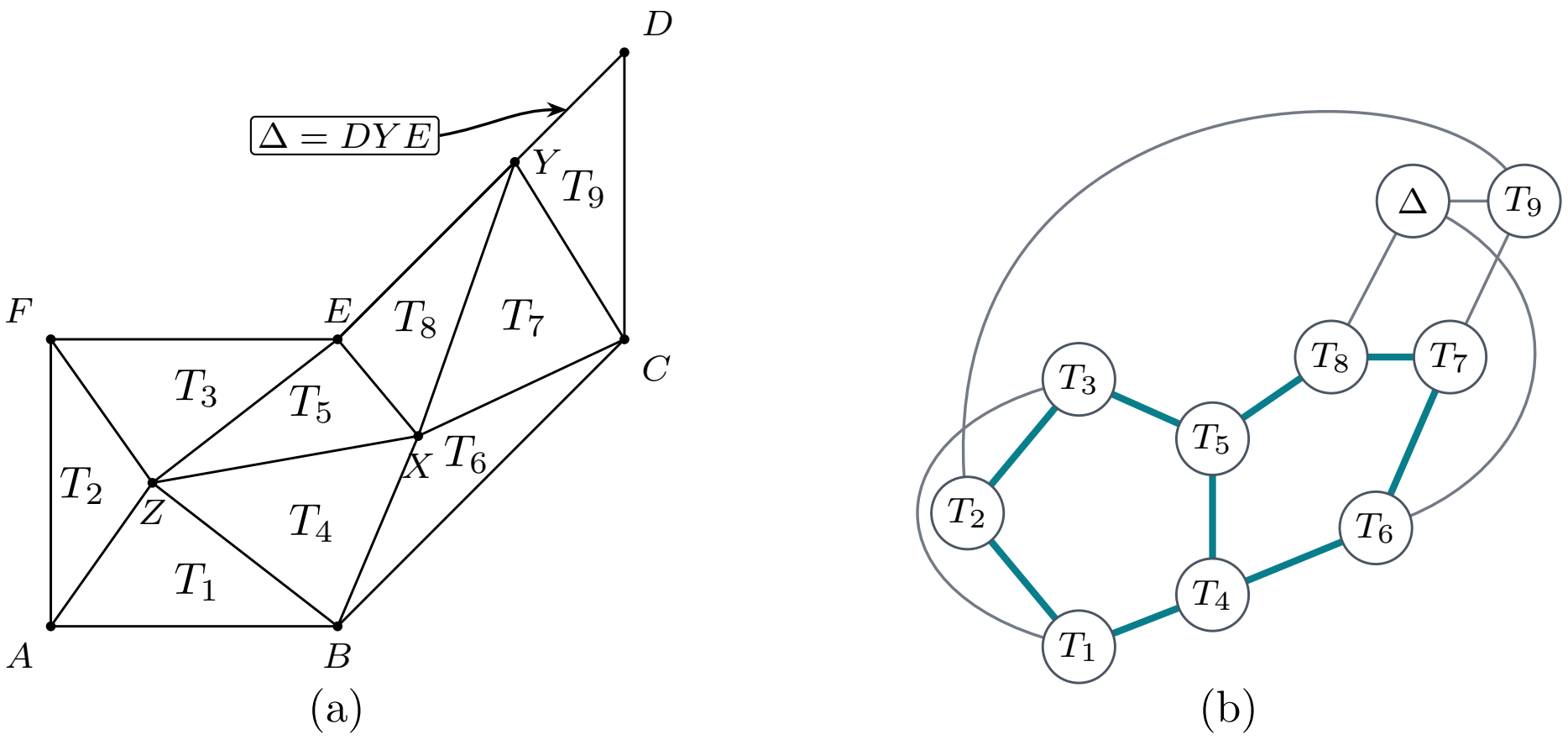
\begin{figure}[!ht]
\centering

\begin{minipage}{\textwidth}
{\centering

\subcaptionbox{\label{fig:balanced-hexagon-flow-a}}{%
\begin{tikzpicture}[x=1.13cm,y=1.13cm,font=\small,
  baseline=(current bounding box.north)]
  \pgfmathsetmacro{\ph}{(1+sqrt(5))/2}

  \coordinate (A) at (0,0);
  \coordinate (B) at (2,0);
  \coordinate (C) at (4,2);
  \coordinate (D) at (4,4);
  \coordinate (E) at (2,2);
  \coordinate (F) at (0,2);
  \coordinate (Z) at ({-9+6*\ph},1);
  \coordinate (X) at ({-12+9*\ph},{-10+7*\ph});
  \coordinate (Y) at ({2*\ph},{2*\ph});

  \draw[black,line width=.55pt]
    (A)--(B)--(C)--(D)--(E)--(F)--cycle
    (A)--(Z)--(B)
    (F)--(Z)--(E)
    (Z)--(X)
    (B)--(X)--(E)
    (C)--(X)--(Y)
    (E)--(Y)
    (C)--(Y);

  \node[draw=black,rounded corners=1pt,fill=white,text=black,
        inner sep=1.5pt,font=\scriptsize] (deglabel) at (2.05,3.42)
        {$\Delta=DYE$};
  \draw[black,-{Stealth[length=1.6mm]},line width=.65pt]
        (deglabel.east) to[out=10,in=180] (3.60,3.60);

  \node at (1.02,.31) {$T_1$};
  \node at (.22,.98) {$T_2$};
  \node at (1.02,1.66) {$T_3$};
  \node at (1.82,.72) {$T_4$};
  \node at (1.82,1.55) {$T_5$};
  \node at (2.90,1.20) {$T_6$};
  \node at (3.30,2.15) {$T_7$};
  \node at (2.55,2.14) {$T_8$};
  \node at (3.72,3.05) {$T_9$};

  \foreach \p/\pos in {A/below left,B/below,C/below right,D/above right,
                       E/above,F/above left,Z/below,X/below,Y/right}
    \fill (\p) circle (1.1pt) node[\pos,font=\scriptsize] {$\p$};
\end{tikzpicture}}
\hspace{1.1cm}
\subcaptionbox{\label{fig:balanced-hexagon-flow-b}}{%
\begin{tikzpicture}[x=1.17cm,y=1.17cm,font=\scriptsize,
  baseline=(current bounding box.north),
  vertex/.style={circle,draw=edgecol,fill=white,line width=.55pt,
                 minimum size=5.7mm,inner sep=0pt},
  ordinary/.style={draw=edgecol!78,line width=.65pt},
  boundary/.style={ordinary},
  marked/.style={draw=markcol,line width=1.55pt}]

  \coordinate (q1) at (.95,.40);
  \coordinate (q2) at (.20,1.30);
  \coordinate (q3) at (.95,2.20);
  \coordinate (q4) at (1.85,.75);
  \coordinate (q5) at (1.85,1.80);
  \coordinate (q6) at (2.95,1.20);
  \coordinate (q7) at (3.45,2.35);
  \coordinate (q8) at (2.65,2.35);
  \coordinate (q9) at (3.95,3.40);
  \coordinate (qd) at (3.20,3.40);

  \draw[marked] (q4)--(q1)--(q2)--(q3)--(q5);
  \draw[marked] (q4)--(q6)--(q7)--(q8)--(q5);
  \draw[marked] (q4)--(q5);

  \draw[ordinary] (q7)--(q9)--(qd)--(q8);

  \draw[boundary]
    (q1) .. controls (-.50,.70) and (-.50,1.90) .. (q3);
  \draw[boundary]
    (q6) .. controls (4.20,1.55) and (4.45,2.90) .. (qd);
  \draw[boundary]
    (q2) .. controls (-.20,4.35) and (3.65,4.45) .. (q9);

  \node[vertex] at (q1) {$T_1$};
  \node[vertex] at (q2) {$T_2$};
  \node[vertex] at (q3) {$T_3$};
  \node[vertex] at (q4) {$T_4$};
  \node[vertex] at (q5) {$T_5$};
  \node[vertex] at (q6) {$T_6$};
  \node[vertex] at (q7) {$T_7$};
  \node[vertex] at (q8) {$T_8$};
  \node[vertex] at (q9) {$T_9$};
  \node[vertex] at (qd) {$\Delta$};
\end{tikzpicture}}
\par}

\begin{center}
\begin{minipage}{.98\linewidth}
\footnotesize
\noindent\textbf{On the left:}
The balanced hexagon $ABCDEF$ is dissected into nine nondegenerate triangles, not of equal area. Adding a tenth degenerate triangle $\Delta=DYE$ allows one to pass to a dual trivalent graph. Let $K=\mathbb Q(\sqrt5)$, $\varphi=(1+\sqrt5)/2$, and $\mathcal{O}_K=\mathbb Z[\varphi]$. All vertices of the dissection are chosen in
$\mathcal O_K^2$; specifically, the vertices of the hexagon are
$A(0,0), B(2,0), C(4,2), D(4,4), E(2,2), F(0,2)$. The remaining vertices are $Z(-9+6\varphi,1), X(-12+9\varphi,-10+7\varphi)$, and $Y(2\varphi,2\varphi)$.

\smallskip
\noindent\textbf{On the right:}
The flow on an edge joining $T_i$ and $T_j$ is the vector of their common side, with opposite boundary sides identified under the pairing given by balancedness. Set $\mathcal O=\mathcal O_{K,(2)}$, take $\pi=2$, and let $P=\mathcal Oe_1\oplus\mathcal Oe_2$ and $L=\mathcal{O}(e_1+e_2)\oplus\mathcal O(2e_2)$, where $e_1:=(1, 0)$ and $e_2:=(0, 1)$. The teal edges form the subgraph $\Gamma_L$; equivalently, they are the edges whose flow values lie in $L\setminus\pi P$. The flow value on every grey edge lies in $\pi P$. For every $v\in \{T_1, T_2, T_3, T_4, T_5, T_7\}$ one has $L_v=L$, whereas for $v\in \{T_6, T_8\}$ the lattice $L_v$ is a strict sublattice of $L$. Thus, $S_L\subsetneq V(\Gamma_L)$. 

\end{minipage}
\end{center}

\end{minipage}

\caption{A dissection of a balanced hexagon and its associated trivalent flow graph.}
\label{fig:balanced-hexagon-flow}
\end{figure}

\begin{proof}
Smith normal form together with the condition $L \nsubseteq \pi P$ implies the existence of a basis $b_1, b_2$ of $P$ such that
\[
P=\mathcal O b_1 \oplus \mathcal O b_2,
\qquad
L=\mathcal O b_1 \oplus \mathcal O \pi^r b_2
\]
for some $r \geq 0$. 

Let us call the edges of $\Gamma_L$ marked edges. We refer to Figure \ref{fig:balanced-hexagon-flow-b} for an example. To see (i), let $v\in V(\Gamma_L)$ and suppose $e$ is a marked edge adjacent to $v$, so that $f(e) \in L \setminus \pi P$. Let $e_2$ and $e_3$ be the other two edges adjacent to $v$. The flow condition rules out that both $f(e_2)$ and $f(e_3)$ lie in $\pi P$, and so we may assume that $f(e_2)\not\in \pi P$. Suppose, for a contradiction, that $f(e_2)\not\in L$. Then in the adapted basis the vector $f(e)$ has unit first coordinate, while $f(e_2)$ has second coordinate of valuation strictly smaller than $r$. Hence $f(e)$ and $f(e_2)$ generate a lattice strictly containing $L$, contrary to the hypothesis. 

Suppose now that a vertex $v$ has degree 2 in $\Gamma_L$. Let $e_1, e_2$ be the two edges of $\Gamma_L$ adjacent to $v$ and let $e_3$ be the edge of $\Gamma$ adjacent to $v$ which is not in $\Gamma_L$. Since $f(e_1), f(e_2)\in L\setminus \pi P$, the flow condition forces $f(e_3)\in L$. Since $e_3$ is not marked, this implies $f(e_3)\in \pi P$, and hence $f(e_3)\in L\cap \pi P$, proving (ii).

We now turn to proving (iii). First, if $v\in S_L$ then the three adjacent edges $e_1, e_2, e_3$ satisfy $f(e_1), f(e_2), f(e_3)\in L$ and, since $L\nsubseteq \pi P$, at least one edge must be marked, for otherwise $L_v$ would be a sublattice of $\pi P$. Thus, $S_L\subseteq V(\Gamma_L)$. Using the basis $b_1,\pi^r b_2$ of $L$, identify $W=L/\pi L$ with $k^2$. Let $B\in(\Lambda^2 k^2)^\vee$ be the residue determinant defined by
\[
B((s_1,s_2),(t_1,t_2))=s_1t_2-s_2t_1,
\]
and let $q(s_1,s_2)=s_1s_2$.
Because the residue field has characteristic $2$,
\[
B(s,t)=q(s)+q(t)+q(s+t).
\]

By (i) and (ii) all three labels at a vertex of $\Gamma_L$ lie in $L$. At a vertex $v$ of $\Gamma_L$, let $w_1,w_2,w_3 \in W$ be the reduced outgoing labels. They satisfy $w_1+w_2+w_3=0$. Under the identification $\det L/\pi\det L\simeq k$ induced by the chosen basis of $L$, the image of the determinant of any two incident labels is
\[
B(w_1,w_2)=q(w_1)+q(w_2)+q(w_3).
\]
If $v \in S_L$, this image is the image of $A_v$. If $v \not\in S_L$, i.e. if $L_v\not = L$, then the incident labels generate a proper submodule of $L$, so the determinant image is zero in $\det L/\pi\det L\simeq k$. Therefore $\sum_{v \in S_L} A_v \bmod \pi \det L$ is the sum, over all vertices of $\Gamma_L$, of the corresponding three $q$-terms.

Every marked edge contributes twice to this sum, once from each endpoint, and therefore cancels in characteristic $2$. At a degree-$2$ vertex the unmarked label lies in $L \cap \pi P$; in the adapted coordinates its first residue coordinate is zero, so its $q$-value is zero. Thus no boundary term remains, and the total sum lies in $\pi \det L$.
\end{proof}

\begin{proposition}[Lifting a residue flow]\label{prop:residue-lift}
Let $\Gamma$ be a graph, let $A$ be an abelian group, and let $A'\subseteq A$ be a subgroup. Suppose that $g:\overrightarrow{E}(\Gamma)\to A$ is an antisymmetric labelling of the oriented edges such that, for every vertex $v$ of $\Gamma$,
\[
\sum_{e\text{ leaving }v} g(e)\in A'.
\]
Then there is an $A$-valued flow\footnote{The notion of a flow extends naturally from trivalent graphs to arbitrary graphs.} $z:\overrightarrow{E}(\Gamma)\to A$ such that
\[
z(e)\equiv g(e)\pmod {A'}
\]
for every oriented edge $e$ of $\Gamma$.
\end{proposition}

\begin{proof}
After reducing modulo $A'$, the labelling $g$ becomes an $(A/A')$-valued flow on $\Gamma$. Thus it is a $1$-cycle with coefficients in $A/A'$. The short exact sequence
\[
0\to A'\to A\to A/A'\to 0
\]
gives an exact sequence
\[
H_1(\Gamma;A)\to H_1(\Gamma;A/A')\to H_0(\Gamma;A')\to H_0(\Gamma;A).
\]
The last map is induced by the inclusion $A'\subseteq A$, hence is injective. Therefore $H_1(\Gamma;A)\to H_1(\Gamma;A/A')$ is surjective. Since a graph has no
$2$-cells, $H_1$ is the group of flows. The reduced flow therefore lifts to an $A$-valued flow $z$, as required.
\end{proof}
\begin{remark}
The homological proof above can be replaced by an elementary spanning-tree argument. Namely, choose a spanning forest of $\Gamma$ and lift the values of the reduced flow on the edges not in the forest in an arbitrary way from $A/A'$ to $A$. The values on the edges of each tree of the forest can then be determined successively, from the leaves toward the roots, by imposing zero divergence. The condition at each root follows automatically by summing the divergences over its connected component.
\end{remark}

We are now ready to prove the main theorem of this section.

\begin{proof}[Proof of Theorem~\ref{thm:flow-main}]
By scaling we may assume that all labels lie in $\mathcal{O}^2$ and then take the ambient lattice to be $P:=\mathcal{O}^2$. Scaling multiplies every $A_v$ by the same nonzero square and hence preserves the assertion. We induct on the value $m_0:=m(\Gamma, f)$, which is non-negative. If $m_0=0$, then every $v$ with $m_v=m_0$ has $L_v=P$, and in particular $L_v\nsubseteq\pi P$; Proposition~\ref{prop:base-case}(iii), applied with $L=P$, gives the result.

We now prove the induction step. Let $m_0=m(\Gamma,f)>0$, and let $\mathcal L$ be the set of lattices $L\subseteq P$ with $\nu_{\mathfrak{p}}(L)=m_0$ and $L\nsubseteq \pi P$. Call an edge primitive if its label lies in $P\setminus \pi P$. Every primitive edge belongs to a unique lattice in $\mathcal L$: if $a\in P\setminus\pi P$, then this lattice is $\mathcal{O}a+\pi^{m_0}P$. Indeed, writing $P=\mathcal Ob_1\oplus\mathcal Ob_2$ and $L=\mathcal Ob_1\oplus\mathcal O\pi^{m_0}b_2$ as in the proof of Proposition~\ref{prop:base-case}, we see that $\pi^{m_0}P\subseteq L$, so $L\supseteq\mathcal Oa+\pi^{m_0}P$; and the two lattices coincide, since both have determinant of valuation $m_0$. Hence the graphs $\Gamma_L$ are edge-disjoint. Also, since $m_0=m(\Gamma, f)=\min_{v\in V} m_v$, no $L_v$ strictly contains a member of $\mathcal{L}$, for otherwise we would have $m_v<m_0$. Thus, by Proposition~\ref{prop:base-case}(i), each $\Gamma_L$ has degree $2$ or $3$ at every one of its vertices, and so two distinct $\Gamma_L$'s cannot meet at a vertex without sharing an edge. Thus the graphs $\Gamma_L$ are vertex-disjoint as well.

Let $T$ be the set of vertices of weight $m_0$ that are incident to a primitive edge. We claim that $T$ is the disjoint union of the sets $S_L$, $L\in\mathcal L$. Indeed, if $v\in T$, then its primitive edge lies in $\Gamma_L$ for a unique $L\in\mathcal L$; by Proposition~\ref{prop:base-case}(i) and (ii) all labels at $v$ lie in $L$, so $L_v\subseteq L$, and $\nu_{\mathfrak p}(L_v)=m_0=\nu_{\mathfrak p}(L)$ forces $L_v=L$, i.e.\ $v\in S_L$. Conversely, every $v\in S_L$ has weight $\nu_{\mathfrak p}(L)=m_0$ and lies in $V(\Gamma_L)$ by Proposition~\ref{prop:base-case}(iii), hence is incident to a primitive edge; and the sets $S_L$ are pairwise disjoint, since $L_v=L$ determines $L$. By Proposition~\ref{prop:base-case}(iii), each sum $\sum_{v\in S_L}A_v$ lies in $\pi\det L$ and so has valuation at least $m_0+1$; therefore $\sum_{v\in T}A_v$ has valuation $>m_0$.

It remains to handle the vertices of weight $m_0$ whose incident labels all lie in $\pi P$. We first clear all primitive edges. Fix $L\in\mathcal L$ and put
\[
M=L\cap \pi P.
\]
On $\Gamma_L$, the restriction of $f$ is an $L$-valued edge labelling whose divergence at every vertex lies in $M$: at a degree-$3$ vertex the divergence is zero, while at a degree-$2$ vertex it is the negative of the outside edge
label, which lies in $M$ by Proposition~\ref{prop:base-case}. By Proposition~\ref{prop:residue-lift}, there is an $L$-valued flow $z_L$ on $\Gamma_L$ such that
\[
z_L(e)\equiv f(e)\pmod M
\]
for every edge of $\Gamma_L$.

Replace $f$ by $f-z_L$ on $\Gamma_L$ and leave it unchanged outside $\Gamma_L$. Since $z_L$ is a flow, the flow condition is preserved. Since the graphs $\Gamma_L$ are vertex-disjoint, this can be done simultaneously for all $L\in\mathcal L$. Afterwards every primitive edge has been moved into $\pi P$, and every non-primitive edge was already in $\pi P$. Thus the new flow takes values in $\pi P$.

Moreover, every vertex touched by some $\Gamma_L$ now has all incident labels in $M=L\cap\pi P$. In an adapted basis
\[
P=\mathcal{O} a\oplus\mathcal{O} b,
\qquad
L=\mathcal{O} a\oplus\mathcal{O} \pi^{m_0}b,
\]
one has
\[
M=\pi\mathcal{O}a\oplus\mathcal{O}\pi^{m_0}b,
\]
so $\nu_{\mathfrak{p}}(\det M)=m_0+1$. Hence every touched vertex has weight at least $m_0+1$ after the clearing.

Now divide the cleared flow by $\pi$. If there is no untouched vertex of weight $m_0$, we are done. Otherwise, all labels at such a vertex lie in $\pi P$, so $m_0\geq 2$. The vertices of minimal weight in the divided flow are exactly the untouched vertices which had weight $m_0$ before division; their new weight is $m_0-2$. By the induction hypothesis, their determinant sum after division has valuation $>m_0-2$. After multiplying back by $\pi^2$, we find that their original contribution has valuation $>m_0$. Together with the primitive contribution already handled, this proves the induction step.
\end{proof}

\section{Proof of Theorem~\ref{thm:main}}\label{sec:proof-of-the-main-theorem}

To deduce Theorem~\ref{thm:main} from Theorem~\ref{thm:flow-main}, we first prove the following simple proposition.

\begin{proposition}\label{prop:algebraic}
If a balanced polygon admits an equidissection into $n$ triangles, then there exists a balanced polygon with an equidissection into $n$ triangles such that all polygon and triangle vertices have real algebraic coordinates.
\end{proposition}

\begin{proof}
Fix such an equidissection of a balanced polygon into $n$ triangles. Its incidences, vertex orders, and non-overlap can be preserved by fixing the signs, including zeros, of all affine determinants of triples of vertices. These conditions, together with equality and positivity of the triangle areas, can be expressed by finitely many polynomial equations and strict inequalities in the vertex coordinates, with coefficients in $\mathbb Q$.
Fix also the original partition of the boundary edges into parallel direction classes. Within each class, parallelism and vanishing of the sum of the oriented boundary vectors are polynomial equations.

Since $\overline{\mathbb Q}\cap\mathbb R$ and $\mathbb R$ are real closed fields, the Tarski--Seidenberg transfer principle (see, e.g., \cite[Proposition~5.2.3]{BochnakCosteRoy1998}) shows that existence of a real solution to this system implies existence of a solution over $\overline{\mathbb Q}\cap\mathbb R$. The resulting realisation is still balanced: each original direction class has zero vector sum, and this remains true if several classes merge. It therefore gives the required balanced polygon and equidissection into $n$ triangles.
\end{proof}

\begin{remark}
Alternatively, one could allow $K$ in Theorem~\ref{thm:flow-main} to be any finitely generated extension of $\mathbb Q$, provided that $\nu_{\mathfrak p}$ is replaced by a suitable valuation. Indeed, the proof of Theorem~\ref{thm:flow-main} uses only that $\mathcal O$ is a discrete valuation ring whose residue field has characteristic~$2$, and such a valuation exists on every finitely generated $K$: choose a transcendence basis $t_1,\dots,t_r$ of $K/\mathbb Q$ and let $\nu$ be any extension to $K$ of the $2$-adic Gauss valuation on $\mathbb Q(t_1,\dots,t_r)$, defined on polynomials by taking the minimum of the $2$-adic valuations of their coefficients. The Gauss valuation is discrete with residue field $\mathbb F_2(\bar t_1,\dots,\bar t_r)$, and since $K$ is a finite extension of $\mathbb Q(t_1,\dots,t_r)$, the valuation $\nu$ is again discrete with residue field of characteristic~$2$. Note that not every extension of the $2$-adic valuation to $K$ is discrete, so the choice of $\nu$ matters. With $\nu$ in place of $\nu_{\mathfrak p}$ and its valuation ring in place of $\mathcal O_{K,\mathfrak p}$, the proof of Theorem~\ref{thm:flow-main} remains unchanged. Taking $K$ to be the field generated by the coordinates of the dissection would then make Proposition~\ref{prop:algebraic} unnecessary.
\end{remark}

\begin{proof}[Proof of Theorem~\ref{thm:main}]
Suppose that a balanced polygon admits an equidissection into an odd number $n$ of triangles. By Proposition~\ref{prop:algebraic} we may assume that all vertices of the polygon and of the triangles lie in $K^2$ for some number field $K$; let $X\subset K^2$ be the set of these vertices, and choose a prime ideal $\mathfrak p\subset\mathcal O_K$ above~$2$.

We first cut the sides of the triangles into pieces that match up. Subdivide the boundary of the polygon at the points of $X$ lying on it. The resulting boundary segments are still balanced, so in each direction the segments oriented one way can be matched against those oriented the other way, subdividing further where necessary; the new subdivision points lie in $K^2$, since the ratio of the lengths of two parallel segments with endpoints in $K^2$ lies in $K$. Let $Y$ be $X$ together with these points, and cut every side of every triangle at the points of $Y$ lying on it. Each piece on the boundary of the polygon is now paired with an opposite piece, and each other piece lies on a side of exactly two triangles.

Finally, whenever a side $ab$ of a triangle has been cut at points $p_1,\dots,p_{k-1}$ (in this order), add the degenerate triangles $ap_1b,\;p_1p_2b,\;\dots,\;p_{k-2}p_{k-1}b$. Every boundary piece is now a side of exactly one genuine or degenerate triangle, and every other side of a genuine or degenerate triangle is a side of exactly two of them.

Let a graph $\Gamma$ have one vertex for each genuine or degenerate triangle, one edge for each side shared by two of them, and one edge for each pair of paired boundary pieces. 
Orient the genuine triangles counterclockwise, and orient each degenerate triangle so that on its longest side it induces the orientation opposite to the one induced by the other triangle having that segment as a side. Then any two triangles sharing a side induce opposite orientations on it, and every boundary piece carries the orientation of the boundary of the polygon, so that paired pieces are opposite vectors. For an edge of $\Gamma$ joining $T$ and $T'$, label its orientation from $T$ to $T'$ by the vector of the corresponding side or boundary piece of $T$, oriented by $T$; the opposite orientation then receives the opposite vector. The three outgoing labels at every graph vertex are therefore the oriented sides of the corresponding triangle and sum to zero. Thus $\Gamma$ is trivalent and these labels define a $K^2$-valued flow. Note that $\Gamma$ is loopless, since a triangle never shares a side with itself and two boundary pieces of one triangle are never paired, but it may have multiple edges (e.g. a parallelogram cut along a diagonal gives two vertices joined by three edges). We refer to Figure \ref{fig:balanced-hexagon-flow} for an example of a dissection and the corresponding trivalent graph. See also \cite[Section~3]{RudenkoTrivalentGraphs} for the details of this construction in the case $K=\mathbb{Q}$.

The resulting graph has vertices of two types: those corresponding to degenerate triangles have $A_v=0$, while for those corresponding to nondegenerate triangles we may choose the sign of $A_v$ so that $A_v=2S$, where $S$ is the common area of the triangles. Thus the latter are precisely the vertices of minimal weight, and their number is $n$. By Theorem~\ref{thm:flow-main},
\[
\nu_{\mathfrak p}(2nS)=\nu_{\mathfrak p}\Bigl(\sum_v A_v\Bigr)>\nu_{\mathfrak p}(2S),
\]
where the sum is over the nondegenerate triangles. Hence $\nu_{\mathfrak p}(n)>0$. Since $\nu_{\mathfrak p}$ restricts to a positive multiple of the $2$-adic valuation on $\mathbb Q$, this means that $n$ is even, a contradiction.
\end{proof}

\bibliographystyle{alpha}

\end{document}